\documentclass[a4paper,12pt]{article}
\usepackage[latin1]{inputenc}
\usepackage[T1]{fontenc}
\usepackage[a4paper,lmargin=4.1cm,textwidth=12.8cm,tmargin=4.9cm,textheight=18.5cm]{geometry}
\usepackage[all,2cell,v2]{xy}
\usepackage[all]{xy}
\UseAllTwocells
\usepackage[square,numbers,sort]{natbib}
\usepackage{ifthen}
\usepackage{mathbbol}
\usepackage{savesym}
\usepackage{amsmath}
\usepackage{amssymb}
\savesymbol{iint}
\restoresymbol{TXF}{iint}
\usepackage[thmmarks,standard]{ntheorem}
\usepackage{microtype}
\usepackage{stmaryrd}
\newboolean{PourEditeur}
\setboolean{PourEditeur}{false}
\ifthenelse{\boolean{PourEditeur}}{%
}{%
  \usepackage{color}
  \definecolor{violet}{rgb}{ .615, .122, 1.0}%
  \usepackage[colorlinks,linkcolor=violet,citecolor=violet,urlcolor=violet]{hyperref}
}
\usepackage{cleveref}

\usepackage[toc,page]{appendix}

\theoremstyle{plain}
\crefformat{result}{result~#2#1#3}
\crefformat{proposition}{proposition~#2#1#3}
\crefformat{remark}{remark~#2#1#3}
\newcommand*{\TC}{\ensuremath{\mathbb{T}\text{-}\mathbb{C}at}}
\newcommand*{\TG}{\ensuremath{\mathbb{T}\text{-}\mathbb{G}\text{r}}}

\newcommand*{\T}{\ensuremath{\mathbb{T}}}

\title{Corrections to the article: Operadic definition of the non-strict cells}
\author{Camell Kachour}
\begin{document}
\maketitle
\vspace*{3.5cm}

\begin{abstract}
In this short notes we propose a new notion of contractibility for coloured $\omega$-operad defined in the article published in Cahiers de Topologie et de G{\'e}om{\'e}trie Diff{\'e}rentielle Cat{\'e}gorique (2011), volume 4. We propose also an other way to build the monad for free contractible coloured $\omega$-operads.
\end{abstract}

\begin{minipage}{118mm}{\small
}\end{minipage}

\hypersetup{%
     linkcolor=blue}%
%
\vspace*{1cm}

\vspace*{1cm}

\section*{Introduction}
 One terminological correction suggested by Steve Lack is to use the more common name \textit{weak higher transformations} instead of \textit{Non-strict cells} which were defined in \cite{kach:nscellsfinal}. 
 
 In \cite{kach:nscellsfinal} we defined a coglobular complex of $\omega$-operads
 
  \[\xymatrix{B^{0}\ar[rr]<+2pt>^{\delta^{1}_{0}}\ar[rr]<-2pt>_{\kappa^{1}_{0}}
  &&B^{1}\ar[rr]<+2pt>^{\delta^{2}_{1}}\ar[rr]<-2pt>_{\kappa^{2}_{1}}
  &&B^{2}\ar@{.>}[r]<+2pt>^{}\ar@{.>}[r]<-2pt>_{}
  &B^{n-1}\ar[rr]<+2pt>^{\delta^{n}_{n-1}}\ar[rr]<-2pt>_{\kappa^{n}_{n-1}}
  &&B^{n}\ar@{.}[r]<+2pt>\ar@{.}[r]<-2pt>&}\]
such that algebras for $B^{0}$ are the weak $\omega$-categories, algebras for $B^{1}$ are the weak
  $\omega$-functors, algebras for $B^{2}$ are the weak $\omega$-natural transformations, etc. However
  Andr{\'e} Joyal has pointed out to us that there are too many coherence cells for each $B^{n}$ when
  $n\geqslant 2$, and gave us a simple example of a natural transformation which cannot be
  an algebra for the $2$-coloured $\omega$-operad $B^{2}$. 
   In this section we propose a notion of contractibilty, slightly different from those used in 
  \cite{bat:monglob,kach:nscellsfinal}. This new approach excludes the counterexample of Andr{\'e} Joyal.
 
Furthermore the main theorem of section 6 in \cite{kach:nscellsfinal} is false. 
I am indebted to Mark Weber for providing us a counterexample. 
However this false theorem has no impact to main ideas of the article \cite{kach:nscellsfinal}. I am indebted 
to Michael Batanin who told us that the techniques of the coproduct of monads could be adapted to substitute technically for the role of this false theorem, and to Steve Lack who gave us the 
precise result and references that we needed for this correction.

{\bf Acknowledgement.} 
I am grateful to Andr{\'e} Joyal and to Mark Weber to have both pointed out to me these 
imperfections.
\section*{Corrections}
Here $\T$ designates the monad of strict $\omega$-categories on $\omega$-graphs. Notions
of $\T$-graphs, $\T$-categories, constant $\omega$-graphs, can be found in \cite{kach:nscellsfinal,lein1:oper}. The category $\TG_{p,c}$ of pointed $T$-graphs over constant 
$\omega$-graphs, and the category $\TC_c$ of $\T$-categories over constant $\omega$-graphs are both defined in \cite{kach:nscellsfinal}.

\begin{definition}
For any $\mathbb{T}$-graph $(C,d,c)$ over a constant $\omega$-graph $G$,
   a pair of cells $(x,y)$ of $C(n)$ has the \textit{the loop property} if: 
   $s^{n}_{0}(x)=s^{n}_{0}(y)=t^{n}_{0}(x)=t^{n}_{0}(y)$
\end{definition}

\begin{remark}
\label{notation-couleurs}
  Suppose $G$ is a constant $\omega$-graph (see section 1.4 of the article \cite{kach:nscellsfinal}). A $p$-cell of $G$ is denoted by $g(p)$ and this notation has the following meaning: The symbol $g$ indicates the "colour", and the symbol $p$ 
   point out that we must see $g(p)$ as a $p$-cell of $G$, because $G$ has to be seen as an $\omega$-graph even though it is just a set.
\end{remark}

\begin{definition}
For any $\mathbb{T}$-graph $(C,d,c)$ over a constant $\omega$-graph $G$,
   we call \textit{the root cells} of $(C,d,c)$, those cells whose arities are the reflexivity of a $0$-cell $g(0)$ of $G$, where here $"g''$ indicates the colour (see section \ref{notation-couleurs}), or in other words, those cells $x\in C(n)$ ($n\geqslant 1$) such that
   $d(x)=1^{0}_{n}(g(0))$.
\end{definition}

Here $1^{0}_{n}$ designates the reflexivity operators of the free strict $\omega$-category $\T(G)$ 
(see also \cite{kach:nscellsfinal}). These notions of \textit{root cells} and \textit{loop condition} 
are the keys for our new approach to contractibility. These observations motivate us to put
  the following definition of what should be a contractible $\T$-graphs $(C,d,c)$. For each integers 
  $k\geqslant 1$, let us note $\tilde{C}(k)=\{(x,y)\in C(k)\times C(k): x\|y$ and $d(x)=d(y)$, 
  and if also $(x,y)$ is a pair of root cells then they also need to verify the 
  \textit{loop property}: $s^{k}_{0}(x)=t^{k}_{0}(y)\}$. Also we put $\tilde{C}(0)=\{(x,x)\in C(0)\times C(0)\}$.
\begin{definition}
 A contraction on the $\mathbb{T}$-graph $(C,d,c)$, is the datum, for all $k\in\mathbb{N}$, of a map
$\tilde{C}(k)\xrightarrow{[,]_{k}}C(k+1)$ such that
\begin{itemize}
\item$s([\alpha,\beta]_{k})=\alpha, t([\alpha,\beta]_{k})=\beta$,
\item$d([\alpha,\beta]_{k})=1_{d(\alpha)=d(\beta)}$.
\end{itemize}
\end{definition}
A $\T$-graph which is equipped with a
contraction will be called contractible and we use the notation
$(C,d,c;([,]_{k})_{k\in\mathbb{N}})$ for a contractible $\T$-graph.
Nothing prevents a contractible $\mathbb{T}$-graph from being equipped with
several contractions. So here $C\TG_c$ is the category of the contractible
$\T$-graphs equipped with a specific contraction, and morphisms of this
category preserves the contractions. One can also refer to the
category $C\TG_{c,G}$, where here contractible $\T$-graphs are only taken over a specific
constant $\infty$-graph $G$. A pointed contractible $\T$-graphs 
 (see section 1.2 of the article \cite{kach:nscellsfinal}) is denoted 
 $(C,d,c;p,([,]_{k})_{k\in\mathbb{N}})$, and morphisms between two pointed contractible 
 $\T$-graphs preserve contractibilities and pointings. The category of pointed contractible 
 $\T$-graphs is denoted by $C\TG_{p,c}$. The categories $\TG_{p,c}$ and $C\TG_{p,c}$ are
 both locally finitely presentable and the forgetful functor $V$ 
  $$\xymatrix{H\dashv V: C\TG_{p,c}\ar[rr]&&\TG_{p,c}}$$
  is monadic, with induced monad $\T_C$ is finitary. 
    
  Also the category $\TC_c$ is
   locally finitely presentable and the forgetful functor $U$
  $$\xymatrix{M\dashv U: \TC_c\ar[rr]&&\TG_{p,c}}$$
  is monadic, with induced monad $\T_M$ is finitary.
   
 A $\T$-category is contractible if its underlying pointed $\T$-graph lies in $C\TG_{p,c}$. Morphisms between
 two contractible $\T$-categories are morphisms of $\T$-categories which preserve contractibilities. Let us write $C\TC_c$ for the category of contractible $\T$-categories. Also consider the pullback in 
 $\mathbb{C}AT$
 
 \[\xymatrix{C\TG_{p,c}\underset{\TG_{p,c}}\times\TC_c\ar[rr]^(.6){p_1}\ar[dd]_{p_2}&&\TC_c\ar[dd]^{U}\\\\
 C\TG_{p,c}\ar[rr]_{V}&&\TG_{p,c}}\]
We have an equivalence of categories
 $$C\TG_{p,c}\underset{\TG_{p,c}}\times\TC_c\simeq C\TC_c$$
 Furthermore we have the general fact (which can be found in the articles \cite{kelly-tansfinite,lack-monadicity}):
 \begin{proposition}[Max Kelly]
Let $K$ be a locally finitely presentable category, and $Mnd_{f}(K)$ the category
of finitary monads on $K$ and strict morphisms of monads. Then $Mnd_{f}(K)$
is itself locally finitely presentable. If $T$ and $S$ are object
of $Mnd_{f}(K)$, then the coproduct $T\coprod S$ is \textit{algebraic}, which means that 
 $K^T\underset{K}\times K^S$ is equal to $K^{T\coprod S}$ and the diagonal of the pullback square 
 \[\xymatrix{K^T\underset{K}\times K^S\ar[rr]^(.6){p_1}\ar[dd]_{p_2}&&K^S\ar[dd]^{U}\\\\
 K^T\ar[rr]_{V}&&K}\] 
is the forgetful functor $K^{T\coprod S}\longrightarrow K$. Furthermore
the projections \\
$p_1: K^T\underset{K}\times K^S\longrightarrow K^S$ and
$p_2: K^T\underset{K}\times K^S\longrightarrow K^T$ 
 are monadic.  
 \end{proposition}
\begin{remark}
According to Steve Lack this result remains true for monads having ranks in the context of
locally presentable category.
\end{remark}
 
 We apply this proposition to the diagram above which shows that 
 $C\TC_c$ is a locally presentable category, and also that the forgetful functor 
  \[\xymatrix{C\TC_c\ar[r]^{O}&\TG_{p,c}}\]
   is monadic. Denote by $F$ the left adjoint of $O$. If we apply the functor $F$ to the coglobular
complex of $\TG_{p,c}$ build in the article \cite{kach:nscellsfinal}

\[\xymatrix{C^{0}\ar[rr]<+2pt>^{\delta^{1}_{0}}\ar[rr]<-2pt>_{\kappa^{1}_{0}}
  &&C^{1}\ar[rr]<+2pt>^{\delta^{2}_{1}}\ar[rr]<-2pt>_{\kappa^{2}_{1}}
  &&C^{2}\ar@{.>}[r]<+2pt>^{}\ar@{.>}[r]<-2pt>_{}
  &C^{n-1}\ar[rr]<+2pt>^{\delta^{n}_{n-1}}\ar[rr]<-2pt>_{\kappa^{n}_{n-1}}
  &&C^{n}\ar@{.}[r]<+2pt>\ar@{.}[r]<-2pt>&}\]
we obtain the coglobular complex of the coloured $\omega$-operads of the weak higher transformations 
with our corrected notion of contractibility

\[\xymatrix{B_C^{0}\ar[rr]<+2pt>^{\delta^{1}_{0}}\ar[rr]<-2pt>_{\kappa^{1}_{0}}
  &&B_C^{1}\ar[rr]<+2pt>^{\delta^{2}_{1}}\ar[rr]<-2pt>_{\kappa^{2}_{1}}
  &&B_C^{2}\ar@{.>}[r]<+2pt>^{}\ar@{.>}[r]<-2pt>_{}
  &B_C^{n-1}\ar[rr]<+2pt>^{\delta^{n}_{n-1}}\ar[rr]<-2pt>_{\kappa^{n}_{n-1}}
  &&B_C^{n}\ar@{.}[r]<+2pt>\ar@{.}[r]<-2pt>&}\]

 \begin{remark}
 It is evident that the $\omega$-operad $B^{0}_{C}$ of Michael Batanin is still initial in the category of  
 contractible $\omega$-operads equipped with a composition system, where our new approach of  contractibility is considered.  
\end{remark}

\bigbreak{}
  \begin{minipage}{1.0\linewidth}
    Camell \textsc{Kachour}\\
    Department of Mathematics,
    Macquarie University\\
    North Ryde,
    NSW 2109
    Australia.\\
    Phone: 00 612 9850 8942\\
    Email:\href{mailto:camell.kachour@gmail.com}{\url{camell.kachour@gmail.com}}
  \end{minipage}

\end{document}